\documentclass[12pt,oneside,english]{amsart}
\usepackage[latin1]{inputenc}
\usepackage{amssymb}

\makeatletter
\newtheorem{theorem}{Theorem}

\newtheorem{conjecture}{Conjecture}

\numberwithin{equation}{section}
\usepackage{babel}

\makeatother
\begin{document}
\baselineskip=17pt

\title[On perfect and near-perfect numbers]{On perfect and near-perfect numbers}

\author{Vladimir Shevelev}
\address{Departments of Mathematics \\Ben-Gurion University of the
 Negev\\Beer-Sheva 84105, Israel. e-mail:shevelev@bgu.ac.il}
\subjclass{11B83}

\begin{abstract}
 We call $n$ a \slshape near-perfect \upshape number, if it is sum of all its proper divisors, except for one of them \enskip(\slshape"redundant divisor"\upshape). It is a special kind of Sierpi\'{n}ski number \cite{1}. We prove an Euclid-like theorem for near-perfect numbers and obtain some other results for them.
\end{abstract}

\maketitle

\section{Introduction }
 A \slshape perfect number \upshape is a positive integer equals to the sum of its proper positive divisors. Denote $\sigma(n)$ the sum of all positive divisors of $n.$ Then $n$ is a perfect number if and only if $\sigma(n)-n=n,$ or
\begin{equation}\label{1.1}
\sigma(n)=2n.
\end{equation}
The smallest perfect number is 6, because 1, 2, and 3 are its proper positive divisors, and 1 + 2 + 3 = 6.
By direct experiment one can obtain some first perfect numbers: 6,\enskip28,\enskip496,\enskip8128,... These experiments, at the first time, were conducted by Euclid. Moreover, he was successful to obtain the following important result.
\begin{theorem} \label{t1}(Euclid)
If $p$ is such prime that also $2^p-1$ is prime, then $n=2^{p-1}(2^p-1)$ is perfect number.
\end{theorem}
It is interesting that 2000 passed before a new large success in the research of perfect numbers. L. Euler was successful to convert the Euclid's theorem for even perfect numbers.
\begin{theorem} \label{t2}(Euler)
Even perfect numbers have the form $n=2^{p-1}(2^p-1),$ where $p$ and $2^p-1$ are primes.
\end{theorem}
\indent Recall that primes of the form $2^p-1$ are called \slshape Mersenne primes.\upshape Up to now it is not known whether exist infinitely many Mersenne primes. Therefore, it is not known whether exist infinitely many even perfect numbers. Not less difficult problem is: whether exists at least one \slshape odd \upshape perfect number? This question is still open as well.\newline
\indent In connection with study of the perfect numbers, it is natural to split all positive integers into three sets: numbers for which $\sigma(n)<2n,$ perfect numbers and numbers for which $\sigma(n)>2n.$ Numbers of the first set are called \slshape deficient, \upshape while numbers of the third set are called \slshape abundant. \upshape \newpage In contrast to the perfect numbers, it is known that both sets of deficient numbers and abundant numbers are infinite. Many other facts and problems on deficient, abundant and perfect numbers one can find in \cite{1}, Chapter B.\newline
\indent A positive integer is called \slshape pseudoperfect, \upshape or \slshape Sierpi\'{n}ski number\slshape \enskip $($cf. \cite{3}$),$ if it is the sum of \slshape some \upshape of its divisors; e.g., $36=1+2+6+9+18.$ In this paper we study Sierpi\'{n}ski numbers of a special kind. We call $n$ a \slshape near-perfect \upshape number, if it is sum of all its proper divisors, except of \slshape one \upshape of them: $d=d(n).$ The latter divisor we call \slshape redundant. \upshape With help of near-perfect numbers, it is suitable, for a given $l,$ to introduce a notion of $(\mathbf{N}\setminus\{l\})$-perfect number. We call $n$ a $(\mathbf{N}\setminus\{l\})$-perfect number if it is perfect number in case of $n$ is not multiple of $l,$  otherwise, if it is near-perfect number with redundant divisor $l.$\newline
The first near-perfect numbers are (cf. our sequence A181595 in \cite{4}):
\begin{equation}\label{1.2}
12, 18, 20, 24, 40, 56, 88, 104, 196, 224, 234, 368, 464, 650, 992,...
\end{equation}
with the redundant divisors (cf. our sequence A181596 in \cite{4})
\begin{equation}\label{1.3}
 4,  3,  2, 12, 10,  8,  4,   2,   7,  56,  78,   8,   2,   2,  32,...
\end{equation}
 \section{Euclid-like theorem for near-perfect numbers}
For a given $k\geq1,$ consider set $\mathcal{P}_k$ of primes of the form: $2^t-2^k-1,$ where $t\geq k+1.$
\begin{theorem} \label{t3}
Number  $n=2^{t-1}(2^t-2^k-1),$ where $2^t-2^k-1\in \mathcal{P}_k,$ is near-perfect number with redundant divisor $2^k.$
\end{theorem}
\bfseries Proof. \mdseries Since $k\leq t-1,$ then $d=2^k$ is a proper divisor of $n.$
 Besides, $\sigma(n)=(2^t-1)(2^t-2^k)$ and, in view of
 $$\sigma(n)-2n=(2^t-1)(2^t-2^k) -2^t(2^t-2^k-1)=2^k,$$
 the theorem follows. $\blacksquare$ \newline
In contrast to the case of perfect numbers, there exist even near-perfect numbers which have not form $n=2^{t-1}p$ with $p\in \mathcal{P}_k.$
 Indeed, consider near-perfect number $n=650$ from (\ref{1.2}). The redundant divisor for it is d(650)=2. Nevertheless, 650 is not expressible in the considered form. But Theorem \ref{t3} makes plausible the following conjecture.
 \begin{conjecture} \label{c1}
 For given $k,$ there exist infinitely many near-perfect numbers with redundant divisor $2^k.$
 \end{conjecture}
 Let
 $$\mathcal{P}=\bigcup _{k=1}^\infty \mathcal{P}_k.$$\newpage
 The first primes from $\mathcal{P}$ are (cf. our sequence A181741 in \cite{4}):
 \begin{equation}\label{2.1}
  3,5,7,11,13,23,29,31,47,59,61,127,191,223,239,...
  \end{equation}
  Note that all Mersenne primes are in the sequence. Indeed, if a prime $p$ has form $p=2^r-1,$ then $p=2^{r+1}-2^r-1.$ On the other hand, it is easy to see that, if $p$ is in the sequence (\ref{2.1}), then the representation $p=2^t-2^k-1,\enskip k>=1,\enskip t>=k+1,$ is unique for it.\newline
  \indent According to Theorem \ref{t3}, near-perfect numbers of the form $n=2^{t-1}(2^t-2^k-1),$ where $2^t-2^k-1\in \mathcal{P}_k,$ we call $\mathcal{P}_k$-\slshape near-perfect, \upshape for a given $k,$ and $\mathcal{P}$-\slshape near-perfect\upshape, if $k$ is not fixed.

  \section{Near-perfect numbers generated by perfect numbers}

  Let us seek near-perfect numbers in the form: $n=2^x m,$ where m is an even perfect number.
  \begin{theorem} \label{t4}
Number $n$ of the form $n=2^x m,$ where $m$ is even perfect number, is near-perfect if and only if either $x=1$ or $x=p,$ where $p$ is prime such that $2^{p-1}$ is the most power of $2$ dividing $m \enskip (2^{p-1}||m).$
\end{theorem}
 \bfseries Proof. \mdseries Since $2^{p-1}||m,$ then, by Theorem \ref{t1}, $m=2^{p-1}(2^p-1),$ where $2^p-1$ is prime.
 Therefore, $n=2^{p+x-1}(2^p-1).$ We have
 $$\sigma(n)-2n=(2^{x+p}-1)2^p-2^{x+p}(2^p-1)=2^p(2^x-1).$$
 This is a proper divisor of $n$ if and only if either $x=p$ or $x=1.\blacksquare$ \newline
  So, every even perfect number $m=2^{p-1}(2^p-1)$ generates two distinct near-perfect numbers $n_1=2m$ and $n_2=2^pm.$ Note that $n_1$ is $\mathcal{P}_{p+1}$- near-perfect, while $n_2$ is not $\mathcal{P}$-near-perfect.\newline
  \indent Quite another type of near-perfect numbers generated by perfect numbers gives the following theorem.
  \begin{theorem} \label{t5}
Number $n$ of the form $n=2^{p-1}(2^p-1)^2,$ where $p$ and $2^p-1$ are prime, is near-perfect.
\end{theorem}
\bfseries Proof.\mdseries We have
$$\sigma(n)-2n=(2^p-1)((2^p-1)^2+(2^p-1)+1)-2^p(2^p-1)^2=2^p-1.$$

Since $2^p-1$ is a proper divisor of $n,$ then $n$ is near-perfect (which generated by perfect number $2^{p-1}(2^p-1)^2$).$\blacksquare$ \newline
\indent Sequence (\ref{1.3}) shows that near-perfect numbers with odd redundant divisors occur very rarely. We conjecture that all near-perfect numbers with odd redundant divisors have form as in Theorem \ref{t5}.\newpage
\begin{conjecture} \label{c2}
 If the redundant divisor for an even near-perfect is odd, then it is Mersenne prime.
\end{conjecture}
\indent Note that, from Theorems \ref{t4}-\ref{t5} it follows that every perfect number is represented as difference of two near-perfect numbers. Indeed, for every perfect number $m,$ numbers $n_2=2^pm$ and $n_3=(2^p-1)m$ are near-perfect, such that $n_2-n_3=m.$\newline
Besides, the following conjecture seems plausible.
\begin{conjecture} \label{c3}
 If number $l$ is not a power of $2,$ then it could be redundant divisor for at most one near-perfect number.
\end{conjecture}
\indent Note that, if Conjectures \ref{c2}-\ref{c3} are true, then we obtain the following Euler-like theorem.
\begin{theorem} \label{t6}
If Conjectures $\ref{c2}-\ref{c3}$ are true, then every even near-perfect number $n$ with odd redundant divisor has form $n=2^{p-1}(2^p-1)^2,$ where $p$ and $2^p-1$ are primes.
\end{theorem}
\bfseries Proof.\mdseries From Conjecture \ref{c2} we conclude that, redundant divisor $d(n)$ has form $d(n)=2^p-1$
with prime $p$ and $2^p-1.$ Now from Conjecture \ref{c3} and Theorem \ref{t5} we conclude that $n$ has the required form. $\blacksquare$ \newline
\bfseries Remark.\mdseries \enskip In 2010 (see sequence A181595 \cite{4}) the author conjectured that all near-perfect numbers are even. In particular, it is easy to show that odd square-free numbers are never near-perfect. However, at the beginning of 2012 Donovan Johnson \cite{2} found an only up to $2\cdot10^{12}$ odd near-perfect number which is $173369889=3^4\cdot7^2\cdot11^2\cdot19^2.$

\;\;\;\;\;\;\;\;

\end{document}